\documentstyle[11pt]{article}

\def\b#1{{\bf #1}}
\def\i#1{{\it #1}}

\begin{document}
\title{Identities for Tribonacci-related sequences}
\author{Mario Catalani\\
Department of Economics, University of Torino\\ Via Po 53, 10124 Torino, Italy\\
mario.catalani@unito.it}
\date{}
\maketitle
\begin{abstract}
\small{We establish some identities relating two sequences that are, as
explained, related to the Tribonacci sequence. One of these sequences bears
the same resemblance to the Tribonacci sequence as the Lucas sequence does
to the Fibonacci sequence. Defining a matrix that we call Tribomatrix,
which extends the Fibonacci matrix, we see that the other sequence is
related to the sum of the determinants of the 2nd order principal minors of
this matrix.}
\end{abstract}

\section{Antefacts}
Let $S_n$ be the generalized Lucas sequence, also called
generalized Tribonacci sequence,
that is
$$S_{n+1}=S_n+S_{n-1}+S_{n-2},\qquad S_0=3,\,S_1=1,\,S_2=3.$$
$S_n$ is sequence A001644 in \cite{sloane}.
Let $\{\alpha,\,\beta,\,\gamma\}$ be the roots of the characteristics
polynomial $x^3-x^2-x-1=0$ (for an explicit expression see \cite{elia}).
Let us assume that $\alpha$ is the real root, $\beta$ and $\gamma$ are
the complex conjugate roots.
We have $\alpha=1.8392286...$, $\vert\beta\vert =\vert\gamma\vert =
0.737353...$ (see \cite{wolfram}).

\noindent
The Binet's formula (see \cite{elia}) is
$$S_n=\alpha^n+\beta^n+\gamma^n,$$
and the ordinary generating function $A(x)$ is
$$A(x)={3-2x-x^2\over 1-x-x^2-x^3}.$$
Let us consider the following matrix, that we might call Tribomatrix,
$$\b{A}=\left [\begin{array}{ccc} 1&1&0\\
1&0&1\\1&0&0\end{array}\right ].$$
The eigenvalues of this matrix are $\{\alpha,\,\beta,\,\gamma\}$.
Using the relationships between eigenvalues and coefficients of the
characteristic equation we have
\begin{enumerate}
\item $$\alpha+\beta+\gamma=1,$$
\item $$\alpha\beta +\alpha\gamma +\beta\gamma =-1,$$
\item $$\alpha\beta\gamma =1.$$
\end{enumerate}
By induction we get
$$\b{A}^n=\left [\begin{array}{ccc}T_{n+1}&T_n&T_{n-1}\\
T_n+T_{n-1}&T_{n-1}+T_{n-2}&T_{n-2}+T_{n-3}\\
T_n&T_{n-1}&T_{n-2}\end{array}\right ],$$
where $T_n$ are the Tribonacci numbers
(sequence A000073 in \cite{sloane})
$$T_n=T_{n-1}+T_{n-2}+T_{n-3},\qquad T_0=0,\,T_1=1,\,T_2=1.$$

\noindent
Then
\begin{eqnarray*}
{\rm tr}(\b{A}^n)&=&S_n \\
&=&T_n+2T_{n-1}+3T_{n-2},
\end{eqnarray*}
where ${\rm tr}(\cdot)$ is the trace operator.

\noindent
From the generating function we get also immediately
$$S_n=
3T_{n+1}-2T_n-T_{n-1}.$$

\noindent
Define
$$C_n=\alpha^n\beta^n +\alpha^n\gamma^n +\beta^n\gamma^n.$$
Then $C_n$ is the sum of the determinants of the principal minors
of order 2 of $\b{A}^n$ and we
obtain
\begin{eqnarray*}
C_n&=&2T_{n+1}T_{n-2}+T_{n+1}T_{n-1}-T_n^2-2T_nT_{n-1}-
T_{n-1}T_{n-3}+T_{n-2}^2\\
&=&-T_n^2+2T_{n-1}^2+3T_{n-2}^2-2T_nT_{n-1}+2T_nT_{n-2}+
4T_{n-1}T_{n-2}.
\end{eqnarray*}
The sequence $C_n$ is sequence A073145 in \cite{sloane}.

\section{A Recurrence for $C_n$}
We have
\begin{eqnarray*}
&&-C_{n-1}-C_{n-2}+C_{n-3}=\\
&=& -\alpha^{n-1}\beta^{n-1}-\alpha^{n-1}\gamma^{n-1}
-\beta^{n-1}\gamma^{n-1}-\alpha^{n-2}\beta^{n-2}-\alpha^{n-2}\gamma^{n-2}\\
&&\quad\quad -\beta^{n-2}\gamma^{n-2}+
\alpha^{n-3}\beta^{n-3}+
\alpha^{n-3}\gamma^{n-3}+\beta^{n-3}\gamma^{n-3}\\&=&
\alpha^{n-3}\beta^{n-3}(1-\alpha\beta-\alpha^2\beta^2)+
\alpha^{n-3}\gamma^{n-3}(1-\alpha\gamma-\alpha^2\gamma^2)\\
&&\quad\quad
+\beta^{n-3}\gamma^{n-3}(1-\beta\gamma-\beta^2\gamma^2).
\end{eqnarray*}
Using relationships among roots we get
\begin{eqnarray*}
1-\alpha\beta-\alpha^2\beta^2&=&\alpha\beta\gamma
-\alpha\beta-\alpha^2\beta^2\\
&=&\alpha\beta(\gamma-1-\alpha\beta)\\
&=&\alpha\beta(\gamma+\alpha\beta+\alpha\gamma+\beta\gamma-\alpha\beta)\\
&=&\alpha\beta(\gamma+\alpha\gamma+\beta\gamma)\\
&=&\alpha\beta\gamma(1+\alpha+\beta)\\
&=&1+1-\gamma\\
&=&2-\gamma.
\end{eqnarray*}
Upon repeating the same calculations for the other quantities we get
\begin{eqnarray*}
-C_{n-1}-C_{n-2}+C_{n-3}&=&2\alpha^{n-3}\beta^{n-3}-
\alpha^{n-3}\beta^{n-3}\gamma
+2\alpha^{n-3}\gamma^{n-3}\\
&&\quad -\alpha^{n-3}\gamma^{n-3}\beta + 2\beta^{n-3}\gamma^{n-3}
-\beta^{n-3}\gamma^{n-3}\alpha\\
&=&2C_{n-3} -\alpha^{n-4}\gamma^{n-4}\alpha\beta\gamma
-\alpha^{n-4}\beta^{n-4}\alpha\beta\gamma\\
&&\quad -\beta^{n-4}\gamma^{n-4}\alpha\beta\gamma\\
&=& 2C_{n-3}-C_{n-4},
\end{eqnarray*}
that is
$$C_{n-1}=-C_{n-2}-C_{n-3}+C_{n-4}.$$
So we got the recurrence
\begin{equation}
C_n=-C_{n-1}-C_{n-2}+C_{n-3},
\end{equation}
with $C_0=3,\,C_1=-1,\;C_2=-1$.

In this way we obtain easily the ordinary generating function for $C_n$
\begin{equation}
A(x)={3+2x+x^2\over 1+x+x^2-x^3}.
\end{equation}

\section{A Recurrence for $C_{2n}$}
\begin{eqnarray*}
C_{2n}&=&-C_{2n-1}-C_{2n-2}+C_{2n-3}\\
&=&
C_{2n-2}+C_{2n-3}-C_{2n-4}+C_{2n-3}+C_{2n-4}\\
&&\quad -C_{2n-5}-C_{2n-4}-C_{2n-5}+C_{2n-6}\\
&=&
C_{2n-2}+2C_{2n-3}-C_{2n-4}-2C_{2n-5}+C_{2n-6}\\
&=&
-C_{2n-2}-3C_{2n-4}+C_{2n-6}+2C_{2n-2}+2C_{2n-3}\\
&&\quad
+2C_{2n-4}-2C_{2n-5}\\
&=&-C_{2n-2}-3C_{2n-4}+C_{2n-6}+2C_{2n-2}\\
&&\quad -2(-C_{2n-3}-C_{2n-4}+C_{2n-5})\\
&=&-C_{2n-2}-3C_{2n-4}+C_{2n-6}+2C_{2n-2}
-2C_{2n-2}\\
&=&-C_{2n-2}-3C_{2n-4}+C_{2n-6}.
\end{eqnarray*}
So we got the recurrence
\begin{equation}
C_{2n}=-C_{2n-2}-3C_{2n-4}+C_{2n-6},
\end{equation}
with $C_0=3,\,C_2=-1,\;C_4=-5$.
The ordinary generating function is
$$A(x)={3+2x+3x^2\over 1+x+3x^2-x^3}.$$

\section{Identities}
Let $n\ge m$. Then
\begin{eqnarray}
S_nS_{n+m}&=&(\alpha^n+\beta^n+\gamma^n)(\alpha^{n+m}+\beta^{n+m}
+\gamma^{n+m})\nonumber\\
&=&\alpha^{2n+m}+\alpha^n\beta^{n+m}+\alpha^n\gamma^{n+m}+ \alpha^{n+m}\beta^n+
\beta^{2n+m}\nonumber\\
&&\qquad +\beta^n\gamma^{n+m}+\alpha^{n+m}\gamma^n+
\beta^{n+m}\gamma^n+\gamma^{2n+m}\nonumber\\
&=&S_{2n+m}+\alpha^n\beta^n(\alpha^m+\beta^m)
+\alpha^n\gamma^n(\alpha^m+\gamma^m)\nonumber\\
&&\qquad +\beta^n\gamma^n(\beta^m+\gamma^m)\nonumber\\
&=&
S_{2n+m}+\alpha^n\beta^n(S_m-\gamma^m)
+\alpha^n\gamma^n(S_m-\beta^m)\nonumber\\
&&\qquad +\beta^n\gamma^n(S_m-\alpha^m)\nonumber\\
&=&S_{2n+m}+S_m(\alpha^n\beta^n +\alpha^n\gamma^n +\beta^n\gamma^n)+
\nonumber\\
&&\qquad -\alpha^m\beta^m\gamma^m
(\alpha^{n-m}\beta^{n-m} +\alpha^{n-m}\gamma^{n-m} +
\beta^{n-m}\gamma^{n-m})\nonumber\\
&=&S_{2n+m}+S_mC_n-C_{n-m}.
\end{eqnarray}
On the other hand let $n<m$. Then everything goes the same with the
exception of the next-to-last line, where we collect in the third sum
$\alpha^n\beta^n\gamma^n$. Then the result is
\begin{equation}
S_nS_{n+m}=S_{2n+m}+S_mC_n-S_{m-n}.
\end{equation}

\section{Consequences}
\begin{enumerate}
\item
If we put $n=n-1$ and $m=1$ then we get
$$S_nS_{n-1}=S_{2n-1}+C_{n-1}-C_{n-2}.$$
\item
If we put $m=n$ we get
$$S_nS_{2n}=S_{3n}+S_nC_n-3.$$
\item
Generally we have
$$S_nS_{nm}=S_{n(m+1)}+S_{n(m-1)}C_n-S_{n(m-2)}.$$
\item
If we put $m=0$ we get
\begin{eqnarray}
\label{eq:quadrato}
S_n^2&=&S_{2n}+3C_n-C_n\nonumber\\
&=& S_{2n}+2C_n.
\end{eqnarray}
\item For the cube we have
\begin{eqnarray*}
S_n^3&=&S_n^2S_n\\
&=&(S_{2n}+2C_n)S_n\\
&=&S_nS_{n+n}+2S_nC_n\\
&=&S_{2n+n}+S_nC_n-C_{n-n}+2S_nC_n\\
&=&S_{3n}+3S_nC_n-3.
\end{eqnarray*}
\item For the 4-th power
\begin{eqnarray}
\label{eq:quarta1}
S_n^4&=&(S_n^2)^2\nonumber\\
&=&(S_{2n}+2C_n)^2\nonumber\\
&=&S_{2n}^2+4C_n^2+4S_{2n}C_n\nonumber\\
&=&S_{4n}+2C_{2n}+4C_n^2+4S_{2n}C_n.
\end{eqnarray}
But also
\begin{eqnarray}
\label{eq:quarta2}
S_n^4&=&S_n^3S_n\nonumber\\
&=&(S_{3n}+3S_nC_n-3)S_n\nonumber\\
&=&S_nS_{n+2n}+3S_n^2C_n-3S_n\nonumber\\
&=&S_{4n}+S_{2n}C_n-S_n+3S_{2n}C_n+6C_n^2-3S_n\nonumber\\
&=&S_{4n}-4S_n+4S_{2n}C_n+6C_n^2.
\end{eqnarray}

Confronting Equation~\ref{eq:quarta1} and Equation~\ref{eq:quarta2}
we get this other identity
\begin{equation}
\label{eq:cienne}
2S_n=C_n^2-C_{2n}.
\end{equation}
\end{enumerate}

\end{document}